\input amstex
\documentstyle{amsppt}
\magnification=1200
\pagewidth{6.5 true in}
\pageheight{9.0 true in}
\NoBlackBoxes
\topmatter
\title
Order ideals and A generalized Krull height theorem
\endtitle
\rightheadtext{Order ideals}
\author
David Eisenbud, Craig Huneke and Bernd Ulrich
\endauthor
\address
Mathematical Sciences Research Institute,
1000 Centennial Dr.,
Berkeley, CA 94720
\endaddress
\email
de\@msri.org
\endemail
\address
Department of Mathematics, University of Kansas,
Lawrence, KS 66045
\endaddress
\email
huneke\@math.ukans.edu
\endemail
\address
Department of Mathematics,
Michigan State University,
East Lansing, MI 48824
\endaddress
\email
ulrich\@math.msu.edu
\endemail
\thanks
All three authors were partially supported by the NSF.
\endthanks

\def\gr{\text{grade}}

\def\R{\Cal R}

\def\ext{\text{Ext}_R}
\def\bight{\text{bight}}
\def\rad{\sqrt}
\def\n{\bold n}
\def\th{{\text{th}}}
\def\rank{{\text{rank}}}
\def\tr{\text{tr}}
\def\Sym{\text{Sym}}

\def\m{\bold m}
\def\ra{\rightarrow}
\def\lra{\longrightarrow}

\def\var{\varphi}



\def\hgt{\text{ht}}
\def\rf{r\!f}

\abstract{Let $N$ be a finitely generated module
over a Noetherian local ring $(R,\m)$.
We give criteria for the height of
the  order ideal $N^*(x)$ of  an element $x\in N$ to be bounded by
the rank of $N$. The
Generalized Principal Ideal Theorem of Bruns, Eisenbud
and Evans says that this inequality  always holds if $x\in \m N$.
We show that the inequality even holds
if the hypothesis becomes true after first extending scalars
to some local domain and then factoring out torsion. We give
other conditions in terms of residual intersections and integral
closures of modules.

We derive information about
order ideals that leads to  bounds on the heights
of trace ideals of modules---even
in circumstances where we do not have the expected  bounds for the heights
of the order ideals!}
\endabstract
\endtopmatter
\document

\centerline{\bf Introduction}
\bigskip

Let $(R,\m)$ be a Noetherian local ring and let $N$ be a finitely generated
$R$-module. For $x\in N$ we define the
\it order ideal \rm of $x$, written   $N^*(x)$, to be the set
of images of $x$ under
homomorphisms $N\to R$.

The classical Krull height theorem (Krull [1928]) says that
$r$ elements of $R$ either
generate an ideal of height at most $r$, or the unit ideal. This may be
interpreted by saying that if $N$ is free of rank $r$, then
the order ideal of any element $x\in \m N$ has height
at most $r$. Eisenbud and Evans [1976] conjectured that
the same statement would be true for any module $N$;
they proved the conjecture for rings containing a field,
and Bruns [1981] subsequently gave a general argument.
This work leaves the question  addressed in
this paper: \medskip

{\narrower \noindent \it Under what circumstances does
the order ideal of a minimal generator of $N$ have height
at most the rank of $N$?\par}

\medskip
\noindent We give criteria to settle this question in many
cases, and use them in turn to prove related results bounding
the heights of some residual intersections and trace ideals. An interesting
feature is the need for studying Rees algebras and integral
dependence of modules.

To get a feeling for the central question, consider the case where $R$ is
a graded ring, $N$ is a graded module that represents a vector
bundle $\Cal E$ on $X=\text{Proj}(R)$, and $N$ is generated
(perhaps only as a sheaf)
by elements of degree 0. In this setting the order ideal of an element
of degree 0 is just the vanishing locus of the corresponding section of the
vector bundle; and there is a section which vanishes in codimension
greater than $r$ if and only if $\Cal E$ admits a sub-bundle isomorphic to
${\Cal O}_X$ if and only if the top Chern class $c_r({\Cal E})$ vanishes
(for information on Chern classes see Fulton [1984, Chapter 3]).

For example, consider the cotangent bundle $\Omega^1_{{\bold P}^{n-1}}$
of projective $n-1$-space. It has no global sections, but
its twist ${\Cal E}=\Omega^1_{{\bold P}^{n-1}}(2)$ is a bundle of
rank $n-1$ generated by
its global sections. The corresponding module $N$
over the polynomial ring $R$ in $n$ variables is the kernel
of the map $R^n(1)\to R(2)$ sending the $i^\th$ generator to
the $i^\th$ variable; it is generated in degree 0.
The exact sequence
$$
0\to \Cal E\to \Cal O^n(1)\to \Cal O(2)\to 0
$$
and the fact that the ``Chern polynomial''
$$
c_t(\Cal E)=1+c_1(\Cal E)t+c_2(\Cal E)t^2+\dots
\in A^*({\bold P}^{n-1}) = \bold Z[t]/(t^n)
$$
is multiplicative show that $c_t(\Cal E)\equiv(1+t)^n/(1+2t)\
\text{mod}(t^n)$.
In particular
$c_{n-1}(\Cal E)=\sum_{i=0}^{n-1}(-1)^i{n\choose n-1-i}2^i$. A little
arithmetic gives
$$
c_{n-1}(\Cal E)=
\cases
0&\text{ if  $n$ is even;}\\ 1& \text{if $n$ is odd.}
\endcases
$$
Thus
we expect $N$ to have an order ideal of height $\geq n$ if
and only if $n$ is even. In Section 2 we present an algebraic
analysis proving this result (and settling a number of related cases.)

A first suggestion of the role of integrality in the theory is
illustrated by a question of Huneke and Koh: They asked
whether the order ideals of elements in the integral
closure of $\m N$ are always bounded by the rank of $N$.
We prove a still more general statement; the following is
a special case of Theorem 3.1:

\proclaim{Theorem} Let $R$ be an affine domain,
let $N$ be a finitely generated $R$-module of rank $r$, and  let $x\in N$.
Let $R\to S$ be a homomorphism from $R$ to a local domain
$(S,\n)$ and write $SN$ for $N\otimes_RS$  modulo torsion.
If the image of $x$ in  $SN$ lies in $\n SN$, then the height of
$N^*(x)$ is at most the rank of $N$.
\endproclaim

In Section 2 we give systematic methods for constructing examples
where $N^*(x)$ has height greater than $\rank(N)$. We illustrate
our methods by constructing, among other things,
 a graded module $N$ of rank
5 over a polynomial ring in 6 variables such that every
homogeneous element of $N$ has order ideal of height at most
5, but $N$ contains inhomogeneous elements with order
ideal of height 6!


To describe the more refined results of the paper we continue with the
assumption that $R$ is local.
A basic construction  of this paper (in a special case)
is that of the \it perpendicular module \rm of $N$:
If $F$ is a free module of minimal rank with a
surjection $\pi:F\to N$ we define $N^\perp =$ coker$(\pi^*)$.
Set $M=N^\perp$.
We observe (Remark 2.2) that the order ideals of generators
of $N$ correspond to the colon ideals of the form
$U:_R M=\{f\in R\mid fN\subset U\} = \text{ann}( M/U)$
where $M/U$ is a cyclic module. Thus the existence of
order ideals of elements of $N$ having  extraordinary height
is the same as the existence of  submodules $U\subset M$
as above with $U:_RM$ of extraordinary height. A
classic argument of McAdam [1983], as generalized
in Section 2, shows that  $M$ is then
integrally dependent on $U$ (see Corollary 1.3),
and thus $\ell(M)<\mu(M)$. Under some circumstances
we show that the condition $\ell(M)=\mu(M)$
is actually necessary and sufficient for all $x\in N$
to have order ideals of height $\leq \rank(N)$ (Proposition 3.6
and Proposition 3.9).

Using the related constructions, we are able to deduce
information about the colon ideals from information on
order ideals and vice versa. An example is the following
special case of Proposition 4.3, which gives conditions
under which all ``residual intersections" of a module
have the right height:

\proclaim{Proposition} Let $R$ be a regular local ring containing
a field and let $M$ be a finitely generated torsion free $R$-module
of rank $e$.
If $s$ is an integer such that
$\ext^{i}(M,R) = 0$ for $2\leq i\leq s-1$
then, for every submodule $U\subsetneq M$ with $\mu(U)-e+1\leq s$,
$$
\hgt(U:_RM)\leq \max\{0,\mu(U)-e+1\}.
$$
\endproclaim

In Section 4 we also investigate order ideals
of modules of low rank.
 We show under mild hypotheses that if $N$ has
rank $\leq 2$, or $N$ is a $k^\th$ syzygy module of
rank $k$, then all elements of $N$
have  order ideals of height $\leq \rank(N)$ (Proposition 4.1);
under somewhat more stringent conditions we get a similar
result for modules of rank 3 (Proposition 4.2).

In the last section we consider the relationship of
trace ideals and order ideals.
In the case of a module of rank 3 (or a $k^\th$ syzygy
module of rank $k+1$) satisfying mild conditions there may
well be order ideals that are too large; but
we prove that the radical of such an order
ideal must contain the whole \it trace ideal \rm of $N$,
defined as the sum of all order ideals
$\tr(N)=\{f(x)\mid \ f\in N^* \text{ and } x\in N \}$ (Proposition 5.4).
Finally we turn to the question of the possible
heights of trace ideals. The surprising result is that we can
give  a stronger bound for the height of the trace ideal
of $N$ if the height of some order ideal exceeds $\rank(N)$
than in the contrary case (Theorem 5.5).

\bigskip

\bigskip
\bigskip\centerline{\bf 1. Rees Algebras\rm}
\bigskip

In this section we recall the general notion of Rees algebra of a module
introduced in our [2000], and provide some results about integral dependence.

Let $R$ be a Noetherian ring and let $M$ be a finitely generated
$R$-module.  If $Q$ is a prime ideal of $R$ we write
$\mu_Q(M)$ for
the minimal number of generators of $M_Q$ over $R_Q$.  When $(R,\m)$ is a
local ring we set $\mu(M) = \mu_{\m}(M)$.  By $-^*$ we denote the functor
$\text{Hom}_R(-,R)$.  We say that an $R$-linear map $f:M \to F$ is a {\it
versal map}
from $M$ to a free module if $F$ is a free $R$-module and $f^*$ is
surjective.  The latter condition means that every $R$-linear from $M$ to a
free $R$-module factors through $f$.  In our [2000, 0.1 and 1.3]
we define the {\it Rees algebra\/} of $M$ to
be $\R(M)=\Sym(M)/(\cap_gL_g)$
where the intersection is taken over all maps $g$ from $M$ to
free $R$-modules, and $L_g$ denotes the kernel of $\Sym(g)$.
 Equivalently, $\Cal R(M)$ is isomorphic to the image of the map
$\Sym(f):\Sym(M) \ra \Sym(F)$, where $f:M \ra F$ is any versal map to a
free module.  The Rees algebra of any finitely generated module exists and
is unique up to canonical isomorphisms of graded $R$-algebras; in fact the
construction is functorial.  On a less trivial note we prove in  our [2000,
1.4 and 1.5] that the above definition gives the usual notion of Rees
algebras of ideals, and that over a ${\Bbb Z}$-torsion free ring $R$ any
embedding $g:M \to G$ into a free module $G$ can be used to define the Rees
algebra $\Cal R(M)$ as the image of $\Sym(g)$.

Let in addition $U \subset L$ be submodules of $M$, let $U', L'$ be the
images of $U,L$ in $\Cal R(M)$, and consider the $R$-subalgebras $R[U']
\subset R[L']$ of $\Cal R(M)$.  According to  our [2000, 2.1] we say $L$ is
{\it integral} over $U$ {\it in} $M$ if the ring extension $R[U'] \subset
R[L']$ is integral; the largest such module $L$ (which exists and is
unique) is called the {\it integral closure} of $U$ {\it in} $M$; finally,
we say $M$ is {\it integral} over $U$ or $U$ is a {\it reduction} of $M$,
if $M$ is integral over $U$ in $M$.  In  our [2000, 2.2] we record the
following ``valuative criterion of integrality'':

\proclaim{Proposition 1.1}  Let $R$ be a Noetherian ring, let $M$ be a
finitely generated $R$-module, let $U \subset L$ be submodules of $M$, and let
$f: M \ra F$ be a versal map from $M$ to a free $R$-module.  The following
are equivalent:

(1) $L$ is integral over $U$ in $M$.

(2)  For every minimal prime $Q$ of $R$, the module $L'$ is integral over
$U'$ in $M'$, where $'$ denotes images in $F/QF$.

(3)  For every map $M \ra G$ to a free $R$-module and for every
homomorphism $R \ra S$ to a domain $S$, the module $L'$ is integral over
$U'$ in $M'$, where $'$ denotes tensoring with $S$ and taking images in $S
\otimes_R G$.

(4)  For every homomorphism $R \ra V$ to a rank one discrete valuation ring
$V$ whose kernel is a minimal prime of $R$, we have $U'= L'$, where $'$
denotes tensoring with $V$ and taking images in $V \otimes_R F$.

(5) For every map $M \ra G$ to a free $F$-module and every homomorphism $R
\ra V$ to a rank one discrete valuation ring $V$, we have $U' = L'$, where
$'$ denotes tensoring with $V$ and taking images in $V \otimes_R G$.
\endproclaim

Let $R$ be a Noetherian local ring with residue field $k$ and let $M$
be a  finitely generated $R$-module.  In  our [2000, 2.3] we define the {\it
analytic spread} $\ell(M)$ of $M$ to be the Krull dimension of $k \otimes_R
{\Cal R}(M)$.  In case $k$ is infinite one has $\ell(M) = \{\mu(U)\big|U$ a
reduction of $M\}$; furthermore $\ell(M) \le \mu(M)$ and equality holds if
and only if $M$ admits no proper reduction.

The next theorem is due to McAdam [1983, 4.1] in the case of ideals.
Various cases with modules are treated in Rees [1987, 2.5], Kleiman-Thorup
[1994, 10.7], Katz [1995, 2.4], and Simis-Ulrich-Vasconcelos [1999, 5.6].  Our
proof is a reduction to the case treated by Rees.  The result plays an
important role in this paper.

\medskip

\proclaim{Theorem 1.2}   Let $R$ be a locally equidimensional universally
catenary Noetherian ring, let $M$ be a finitely generated $R$-module and
let $U$ be a submodule of $M$ generated by $t$ elements.  If there exists a
minimal prime $Q$ of $R$ such that $M/QM$ is not integral over the image of
$U$ in $M/QM$, then
$$\hgt(U :_R M) \le \max \{0,t + 1-\mu_Q(M)\}.$$
\endproclaim

\demo{Proof} We localize to assume that $R$ is local and equidimensional.
Write $R' = R/Q, M' = M/QM$, and let $U'$ be the image of $U$ in $M'$.  We
have $\hgt(U :_R M) \le \hgt(U' :_{R'} M')$, because $R'(U :_R M) \subset
U' :_{R'} M'$ and $R$ is an equidimensional catenary local ring.  Thus we
may replace $R,U,M$ by $R',U',M'$ to assume that $R$ is a universally
catenary local domain and $M$ is torsion free.  We may assume that $U :_R M
\neq 0$, and hence that rank $U =$ rank $M$.  With these assumptions
the assertion was proved by Rees [1987, 2.5].
\qed
\enddemo

\proclaim{Corollary 1.3} Let $R$ be a local equidimensional universally
catenary Noetherian ring, let $M$ be a finitely generated $R$-module and
let $U$ be a proper submodule of $M$ generated by $t$ elements.  If
$\ell(M) = \mu(M),$ then there exists a minimal prime $Q$ of $R$ such that
$M/QM$ is not integral over the image of $U$ in $M/QM$.  For any such $Q$,
$$ \hgt(U :_R M) \le \max \{ 0,t + 1 - \mu_Q(M)\}.$$
\endproclaim

\demo{Proof}  We may assume that the residue field of $R$ is infinite.  It
follows from our hypothesis that $M$ is not integral over $U$.  By
Proposition 1.1 and the functoriality of the Rees algebra, there exists a
minimal prime ideal $Q$ such that $M/QM$ is not integral over the image of
$U$.  The assertion now follows from Theorem 1.2.
\qed
\enddemo

\bigskip\centerline{\bf 2. Perpendicular Modules\rm}
\bigskip

Central to this paper is the following:

\remark{\bf Definition 2.1} Let $R$ be a Noetherian ring and let $N$ be a
finitely generated
$R$-module with a choice of generators $x_1, \dots ,x_n$. Map a free
$R$-module
with basis $e_1, \dots, e_n$  to $N$ by sending $e_i$ to $x_i$ and denote
this map by $\pi$. We define
$[x_1, \dots ,x_n]^{\bot} = \text{Coker}(\pi^*)$, and write $x_i^{\bot}$
for the
image of $e_i^*$ in $\text{Coker}(\pi^*)$.
When $R$ is local  and the $x_1, \dots , x_n$ are minimal generators, we set $
N^{\bot} = [x_1, \dots ,x_n]^{\bot}$ and call it the
\it perpendicular module \rm to $N$ (indeed, this module only depends on
$N$).
\endremark

\bigskip

\remark{\bf Remark 2.2. Perpendicular modules and colons}
With notation as in 2.1, let $M=[x_1, \dots ,x_n]^\perp$
and let $U$ be the submodule of $M$ generated
by  $x_1^{\perp}, \dots, x_{i-1}^{\perp}, x^{\perp}_{i+1}, \dots,
x^{\perp}_n$. We have
$$
N^*(x_i)=U:_R M.
$$
The right hand side is clearly equal to the $i^\th$ row ideal
(the ideal generated by the elements of the $i^\th$ row)
of any matrix  presenting $M$ with respect to the generating
set $x_1^\perp, \dots, x_n^\perp$.
The $i^\th$ row ideal  equals  $N^*(x_i)$ because the image of
an element $f$ of $N^*$ under $\pi^*$ has $i^\th$ component equal to
$f(x_i)$.
\endremark
\smallskip

An immediate consequence is that if $R$ is local and $N$ has no
free summand, then $\mu(M)=n$. More generally, if $Q$ is any prime
of $R$  we have $\mu_Q(M)=n-\rf_Q(N)$, where $rf_Q(M)$ denotes the
maximal rank of an $R_Q$-free direct summand of $M_Q$.
\medskip

Combining the above remark with Theorem 1.2 we obtain the following
consequence
for order ideals.

\proclaim{Corollary 2.3} Let $R$ be a locally equidimensional universally
catenary Noetherian ring, let $N$ be a finitely generated $R$-module, and
let $Q$ be a minimal prime of $R$.  Choose a generating set $x_1, \dots ,
x_n$ of $N$,
write $M = [x_1, \dots ,x_n]^{\bot},$ and let
$U = Rx_1^{\bot} + \cdots + Rx^{\bot}_{n-1} \subset M$.  If $M/QM$ is not
integral over the image of $U$, then  $\hgt(N^*(x_n)) \leq  r\!f_Q(N).$
\endproclaim

\demo{Proof} Remark 2.2 shows that $N^*(x_n) = U:_R M$, and by Theorem 1.2,
$\hgt(U:_R M) \leq n-\mu_Q(M)$.  Finally, $n- \mu_Q (M) = r \! f_Q(N)$.
\qed
\enddemo

The following version of the
semicontinuity theorem for heights of ideals in a family will
be useful in bounding the heights of order ideals:

\proclaim{Proposition 2.4}
Let $(R,\m)$ be an equidimensional universally catenary
Noetherian local ring, let $\underline{Z} = Z_1,...,Z_n$ be indeterminates
over
$R$, and write $R' = R(\underline{Z})$.
\roster
\item If $\var$ is a matrix over $R$ with $n$ rows, then
$\hgt((\underline{a}\cdot \var)R)\leq
\hgt((\underline{Z}\cdot \var)R')$ for every vector $\underline{a}$ of $n$
elements in $R$.
\item If $N$ is an $R$-module with generating set $x_1, \dots ,x_n$ and
$y = \sum_{i=1}^nZ_ix_i\in N\otimes_R{R'}$, then $\hgt(N^*(x))\leq
\hgt((N\otimes_R{R'})^*(y))$ for every $x\in N$.
\endroster
\endproclaim

\demo{Proof} To see (1) write
$S = R[\underline{Z}]_{(\m, \{Z_i-a_i\})}$. The ring
$R/(\underline{a}\cdot \var)$ is obtained from
$S/(\underline{Z}\cdot \var)$
by factoring out the ideal generated by the $S$-regular sequence $Z_i-a_i,
1\leq i\leq n$.
Since $S$ is equidimensional and catenary, it follows that
$\hgt((\underline{a}\cdot \var)R)\leq \hgt((\underline{Z}\cdot \var)S)$.
Localizing further we deduce (1).

To prove (2) we apply (1) to a matrix $\var$ presenting
$[x_1, \dots ,x_n]^{\bot}$ with
respect to the generating set $x_1^{\bot}, \dots ,x_n^{\bot}$, and invoke
Remark 2.2.
\qed
\enddemo

Let  $N$ be a finitely generated module
over a Noetherian ring $R$.  We say that $N$ satisfies $G_s$, where $s$ is
a positive integer, if $N_Q$
is free of constant rank $r$ for every minimal prime $Q$ of $R$ and
$\mu_Q(N) \le \dim(R_Q) + r-1$ for every prime $Q$ of $R$ with $1 \le
\dim(R_Q) \le s-1$. In case $G_s$ holds for every $s$, the module $N$ is
said to satisfy $G_{\infty}$.

\medskip

We say that $N$ has a \it rank \rm
and write rank$(N) = r$ if $N_Q$ is free of rank $r$ for every associated
prime $Q$ of $R$. The module $N$ is said to be \it orientable \rm if
$N$ has rank $r$ and $(\Lambda^r N)^{**}\cong R$. This differs slightly
from the definition given in Bruns [1987]. Our definition
implies that
$N_Q$ is the direct sum of a free $R_Q$-module and a torsion module for
every prime $Q$ such that depth$(R_Q)\leq 1$. Thus if two modules in a
short exact sequence are orientable, then so is
the third as long as the right-hand module is torsion free locally in depth
one.

\proclaim{Proposition 2.5}
With notation as in 2.1, let $M=[x_1, \dots ,x_n]^\perp$. The module
$N$ has a rank or is orientable if and only if $M$ has
the same property.
\endproclaim

\demo{Proof} We use the exact sequence $$ 0\ra
N^*\overset{\pi^*}\to\lra (R^n)^*\lra M\to 0\tag{*} $$ from which it
follows that $M$ has a rank if and only if $N^*$ has a rank. If $Q$ is
an associated prime of $R$, then $N_Q$ is free if and only if $N^*_Q$
is free. (To see this, notice that if $f\in N^*_Q$ is a free
generator, then the image of $f:\, N_Q\to R_Q$ is faithful, and hence
$f$ is surjective.)  Thus $N^*$ has a rank if and only if $N$ does.

Since $M$ is torsionfree, (*) shows that $M$ is orientable if and only if
$N^*$ is. But $N^*$ is orientable if and only if $N$ is.
\qed
\enddemo

Finally, we remark that every perpendicular module is torsionless
(contained in a free module). Conversely, any finitely generated
torsionless module is the perpendicular module of a finitely
generated torsionless module with respect to some set of generators:
If the torsionless module $N$ is generated by
$x_1,\dots,x_n$ then $N=[x_1^\perp,\dots,x_n^\perp]^\perp$.
(In general, $[x_1^\perp,\dots,x_n^\perp]^\perp$ is the image of
the natural map $N\to N^{**}$.)


\bigskip\bigskip
\centerline{\bf 3. Principal Ideal Theorems\rm}

\bigskip

We are now ready to prove our first main result.  Recall that $rf_Q(N)$
denotes the maximal rank of a free summand of $N_Q$.

\proclaim{Theorem 3.1. Generalized Height Theorem} Let $R$ be a locally
equidimensional universally catenary Noetherian ring, let $N$ be a
finitely generated $R$-module, and let $x \in N$.  Let $R \ra S$ be a
homomorphism of rings
from $R$ to some Noetherian local domain $(S,\n)$ and write $SN$ for
$N\otimes_R S$ modulo $S$-torsion.  If the image of $x$ in $SN$ lies in $\n
SN$, then
$$\hgt(N^*(x)) \leq  \min_Q\{\rf_Q(N)\}$$
where the minimum is taken over all minimal primes $Q$ of $R$ mapping to
zero in $S$.
\endproclaim

\demo{Proof} Replacing $S$ by a rank one discrete valuation ring $V$
containing $S$ and centered on $\n$ we may assume that $S=V$.  Notice that the
order ideal of $x \otimes 1$ in the $V$-module $N \otimes_R V$ is a proper
ideal.

Choose a generating set $x_1, \dots,x_n\!=\!x$ of $N$ and a presentation $R^m
\overset{\psi}\to\lra R^n$ of $N$ with  \linebreak respect to this
generating set.  Consider the perpendicular module $M = \text{Im}(\psi^*) =
[x_1, \dots ,x_n]^{\bot} \subset R^{m*}$ and its submodule $U = Rx^{\bot}_1 +
\cdots + Rx^{\bot}_{n-1}$.  Let $Q$  be a minimal prime of $R$ mapping to
zero in $S = V$.  According to Corollary 2.3 it suffices to prove
that  $M/QM$ is not integral over the image of $U$.

In fact, writing $M' = \text{Im}(\text{Hom}_R(\psi,V)) = \text{Im}(\psi^*
\otimes_R V) \subset R^{m^*} \otimes_R V$, we have $M' = [x_1 \otimes 1,
\dots, x_n \otimes 1]^{\bot}$ and $M/QM$ maps to $M'$ with the image of
$x^{\bot}_i$ being sent to $(x_i \otimes 1)^{\bot}$.  Thus by the
functoriality of Rees algebras it suffices to show that $M'$ is not
integral over $U' = V(x_1 \otimes 1)^{\bot} + \cdots + V(x_{n-1} \otimes
1)^{\bot}$, or equivalently that $U' \neq M'$ (see for instance
Rees [1987, 2.5] or Proposition 1.1).  However according to Remark 2.2,
$U' :_V M'$ is the order ideal of $x \otimes 1$ in $N \otimes_R V$.  Since
this
ideal is proper we conclude that $U' \neq M'$.
\qed
\enddemo

Theorem 3.1 says that whenever the height of the order
ideal of $x$ exceeds the expected value, then the injection of $Rx$ into
$N$ must
be `valuatively split', meaning that after passing to an arbitrary
valuation, the induced map does split. Put differently, either the order
ideal of $x$ has the expected height or else the height of the order ideal
over any valuation ring becomes infinite.  The set of elements having this
property,
but not splitting themselves, is a remarkable class, and exactly the
class we wish to study.

\proclaim{Corollary 3.2}  Let $R$ be a locally equidimensional universally
catenary
Noetherian ring and let $N$ be a finitely
generated $R$-module.  Let $x \in N$ and suppose that $\hgt(N^*(x)) >
\max_Q \{r\!f_Q(N)\}$, where the maximum is taken over all minimal primes
$Q$ of $R$.
\roster
\item   If $N' \ra N$ is an arbitrary epimorphism of $R$-modules and $x'
\in N'$ is an element mapping to $x$, then for every integer $j$,
$\text{Fitt}_{j+1}(N')$ is integral over $\text{Fitt}_j(N'/Rx')$.

\item If $I$ is an arbitrary integrally closed ideal in $R$, then the
natural map  \linebreak $(Rx)\otimes_R(R/I)\ra N\otimes_R(R/I)$ is injective.
\endroster
\endproclaim

\demo{Proof}  To prove (1), we need to verify that for every
homomorphism from
$R$ to a discrete valuation ring $V$,  $\text{Fitt}_{j+1} (N')V =
\text{Fitt}_j(N'/Rx')V$, or equivalently,
$\text{Fitt}_{j+1} ( N'\otimes_R V) = \text{Fitt}_j(N'\otimes_R V/V(x' \otimes
1))$.  However by Theorem 3.1, the image of $x$ generates a free
$V$-summand of rank one in $VN$. Thus  $x'\otimes 1$ generates a free
$V$-summand of rank one in $N'\otimes_R V$, and the asserted equality of
Fitting ideals is obvious.

To prove (2) we have to show that $IN\cap Rx \subset Ix$. Let
$x_1, \dots ,x_n\!=\!x$
be  a generating  \linebreak set of $N$. If $r\in R$ and $rx\in IN$, then
$rx = \sum_{i=1}^n s_ix_i$ for some $s_i\in I$.  Hence $s_1 x_1 + \cdots +
s_{n-1} x_{n-1} + (s_n - r)x = 0$.  Take $N'$ to be the $R$-module
presented by the transpose of the vector $[s_1, \dots ,s_{n-1},s_n-r]$, and
let $x'$ be the last generator.    Part (1) proves that $s_n-r$ is in the
integral
closure of the ideal $(s_1, \dots ,s_{n-1})$. Therefore $r\in I$, since
$s_i\in I$ and $I$ is integrally closed. Hence $rx\in Ix$ as asserted.
\qed\enddemo

\noindent{\bf Remark 3.3.}  In terms of matrices,
Corollary 3.2(1) can be stated as
follows:  Let $\underline{x} = x_1, \dots ,x_n$ be a generating set of $N$
with
$x_n = x$, let $\psi$ be a matrix with $n$ rows satisfying
$\underline{x}\cdot\psi = 0$,
and let $\psi'$ be the matrix obtained from $\psi$ by deleting the last row.
Then for every integer $i$, $I_i(\psi)$ is integral over $I_i(\psi')$.

\medskip

The second corollary answers in the affirmative a question asked
by the second author and Jee Koh in the late 1980's.

\proclaim{Corollary 3.4} Let $(R,\m)$ be an equidimensional universally
catenary Noetherian local ring and let $N$ be a finitely
generated $R$-module.  If $x$ lies in $\overline{\m N},$ the integral
closure of $\m N$ in $N$,  then $\hgt(N^*(x))\leq \min \{\mu_Q(N)\}$, where
the minimum is taken over all minimal primes $Q$ of $R$ so that $N_Q$ is free.
\endproclaim

\demo{Proof} Let $Q$ be a minimal prime of $R$ so that $N_Q$ is free.
Choose a local embedding $R/Q \hookrightarrow V$, where $(V,\n)$ is a rank
one discrete valuation ring, and a versal map $f:N \ra F$ from $N$ to a
free module.  As $x$ lies in the integral closure of $\m N$ in $N$,
Proposition 1.1 shows that $x' \in \m N'$, where $'$ denotes tensoring with
$V$ and taking images in $F\otimes_RV$.  On the other hand since $N_Q$ is
free,
$f_Q$ is split injective and therefore $N' \cong VN$ as defined in Theorem
3.1.  Thus the image of $x$ lies in $\m VN \subset \n VN$, and Theorem 3.1
immediately gives the conclusion.
\qed\enddemo

Of course, Theorem 3.1 gives as an immediate corollary the theorem of
Bruns, Eisenbud, and Evans
in case the ring is equidimensional and catenary. The usual proofs reduce
to this case, and from the paper of Bruns one obtains that if $(R,\m)$ is a
local Noetherian ring,
$N$ is a finitely generated $R$-module, and  $x\in \m N$, then
$$
\hgt(N^*(x)) \leq \max_Q   \{\mu_Q(N)\}
$$
where the maximum ranges over the minimal primes $Q$ of $R$. However, we can
observe that the proof reduces at once to the complete domain case and
proves even a stronger result:

\proclaim {Theorem 3.5} Let $(R,\m)$ be a Noetherian local ring and let
$N$ be a finitely generated $R$-module. If $x\in mN$, then
$\dim(R/N^*(x)) \geq \max_Q \{\dim (R/Q)- \rf_Q(N)\}$
where the maximum is
taken over all primes ideals in $R$.
Thus,
$$
\hgt(N^*(x)) \leq \min_Q\{\dim (R) -\dim(R/Q) + \rf_Q(N)\}.
$$
\endproclaim

\demo{Proof} If $R$ is a complete local domain and $x\in \m N$,
then $\hgt(N^*(x))\leq \text{rank}(N)$ by Theorem 3.1. Now the second
assertion follows, and implies the first statement at
once.

To prove the theorem for any Noetherian local ring, we only need to show
the first inequality.  We first reduce to the case where $R$ is complete by
passing to the completion $\hat R$ of $R$.
Let $Q$ be a prime in $R$ and choose a minimal prime $\hat Q$ over $Q\hat
R$ such
that dim$(\hat R/\hat Q) = \text{dim}(R/Q)$. Since the map $R_Q \to
\hat{R}_{\hat{Q}}$ is flat and local, it
follows that $\rf_Q(N)=  \rf_{\hat Q}(\hat N)$. If the result holds for
$\hat R$ and $\hat N$,  we obtain
$$
\dim(R/N^*(x)) = \dim(\hat R/\hat N^*(x))\geq \text{dim}(\hat R/\hat Q)
- \rf_{\hat Q}(\hat N)= \text{dim}(R/Q) - \rf_Q(N).
$$
Henceforth we assume that $R$ is complete.

Choose an arbitrary prime ideal $Q$ in $R$ and let $J \subset R$ be the
preimage of $(N/QN)^*(x)$,  the order ideal of the image of $x$ in
the $R/Q$-module $N/QN$. Clearly $N^*(x) \subset J$.   Now $\dim(R/N^*(x))
\ge \dim R/J \ge \dim(R/Q) - \rf_Q  (N)$, where the last inequality follows
from the complete
domain case.
\qed\enddemo

Next we wish to find conditions on a module $N$ over a local ring
which guarantee that $N^*(x)$ has the expected height for {\it  every}
$x\in N$.  Obviously, such a module should not have any nontrivial free
summands,
which means that $\mu(N^{\bot}) = \mu(N)$.
We can turn this necessary condition into a sufficient one if we replace
$\mu(N^{\bot})$ by $\ell(N^{\bot})$:

\proclaim {Proposition 3.6} Let $R$ be an equidimensional universally
catenary Noetherian  local ring
and let $N$ be a finitely generated $R$-module.  If
$\ell(N^\perp) = \mu(N)$, then for every $x \in
N$,
$$ \hgt(N^{\ast}(x)) \leq \max_Q \{\rf_Q (N)\}. $$
where the maximum is taken over all minimal primes $Q$ of $R$.
\endproclaim

\demo{Proof} By Theorem  3.1 we may assume that $x$ can be extended to
a minimal generating set $x_1, \dots ,x_n\! =\! x$ of $N$. Write
$U = Rx_1^\perp+\cdots+Rx_{n-1}^\perp
 \subset M = [x_1, \dots ,x_n]^{\bot}\cong N^{\bot}$.
As $\ell(M) = n$, $M$ cannot be integral over $U$.  Hence by Proposition
1.1 and the functoriality of the Rees algebra, there exists a minimal prime
$Q$ of $R$ such that $M/QM$ is not integral over the image of $U$.  The
assertion now follows from Corollary 2.3.
 \qed\enddemo

Here is a result bounding the analytic spread of a module
from below. We will ultimately use it to give a partial converse of
Proposition 3.6.

\proclaim {Proposition 3.7} Let $R$ be a Noetherian local ring with
infinite residue field and let $M$ be a finitely generated torsion free
$R$-module such that $M_Q$ is free of rank $e$ for every minimal prime $Q$
of $R$. Let $U$ be the submodule of $M$ generated
by $t$ general linear combinations of a set of generators of $M$.
If $\ell(M)\leq t$ and $M$ satisfies $G_{t-e+2}$, then
$ht(U:M)\geq t-e+2$.
\endproclaim

Proposition 3.7 follows at once from the following more general
version, which we phrase as a lower bound for the analytic spread
of a module.

\proclaim {Proposition 3.7bis} Let $R$ be a Noetherian local ring with
infinite residue field and let $M$ be a finitely generated torsionless
$R$-module.
Let $X \subset Spec(R)$ be the nonfree locus of
$M$, and assume that for integers $e$ and $t$
we have  $\mu_{P}(M) \leq \dim (R_{P})+ e - 1$
whenever $P \in X$ and
$\dim(R_P)\leq t-e+1$.
Let $U$ be the submodule of $M$ generated
by $t$ general linear combinations of a set of generators of $M$.
If $\hgt(U:_RM)\leq t-e+1$, then
$\ell(M) \geq t+1$.
\endproclaim

\demo{Proof}  Let $x_{1}, \dots , x_{t}$ be the general
elements of $M$ that generate $U$. Using basic element theory and induction
on $i,~ 0 \le i \le t,$ one can show that
$\mu_{P}(M/R x_{1}+ \cdots + Rx_{i})
\leq \dim (R_{P}) + e - i - 1$
for every $P \in X$ with  $i-e+1 \leq
\dim(R_P) \leq t-e+1$.
In particular $M_{P} = U_{P}$ for every $P \in X$ with
$\dim (R_{P}) = t - e +1$. Now suppose that $\ell(M) \leq t$. Then $U$ is a
reduction of $M$ and hence $M_{P} = U_{P}$ for every prime $P \not\in X$.
Thus
$\hgt(U:_RM) > t-e+1$, which yields a contradiction.
\qed\enddemo

The following consequence of Proposition 3.7bis shows that
under good circumstances the analytic spread
is monotonic for inclusions:

\proclaim{Corollary 3.8} Let $R$ be an
equidimensional universally catenary Noetherian local ring
and let $M$ be a finitely generated $R$-module such that $M_Q$ is free of
rank $e$ for every minimal prime $Q$ of $R$.  Write $\ell = \ell(M)$ and
assume that $M$ satisfies $G_{\ell-e+1}$.  If $M'\subset M$ is any
submodule with $\text{codim}(M/M')\geq \ell-e+1$, then
$\ell(M')\geq\ell(M)$.
\endproclaim

\demo{Proof} We may assume that the
residue field of $R$ is infinite.
Let $t=\ell-1$, and let $U$  be the submodule
generated by $t$ general linear combinations of a
set of generators of $M'$. Since $t<\ell(M)$, the
module $M$ is not integral over $U$. Thus by Proposition 1.2, there exists
a minimal prime $Q$ of $R$ such that
$\hgt(U :_R M)\leq \max\{0, t+1-\mu_Q(M)\}\leq t-e+1$.
Since $\text{codim}(M/M')\geq t-e+2$ we have
$\hgt(U :_R M')= \hgt(U :_R M)\leq t-e+1$. Now Proposition 3.7bis gives
$\ell(M')\geq t+1=\ell(M)$.
\qed\enddemo

Here is the promised partial converse of Proposition 3.6:

\proclaim {Proposition 3.9} Let $R$ be a Noetherian local ring with
infinite residue field and let $N$ be a finitely generated $R$-module
such that $N_Q$ is free of rank $r$ for every minimal prime $Q$ of $R$.
Assume that  $N^{\bot}$ satisfies $G_{\infty}$.
If $\hgt(N^*(x))\leq r$ for every $x\in N$, then $\ell(N^\perp) = \mu(N)$.
\endproclaim

Proposition 3.9 is an immediate consequence of the following more general
result.

\proclaim{Proposition 3.9bis} Let $R$ be a Noetherian local ring with
infinite residue field and let $N$ be a finitely generated $R$-module.  Let
$X \subset \text{Spec}(R)$ be the nonfree locus of $N$, and assume that for
an integer $r$ we have $\rf_P(N) \ge r + 1 - \dim(R_P)$ whenever $P \in X$.
 If $\hgt(N^*(x)) \le r$ for every $x \in N$, then $\ell(N^{\perp}) = \mu(N)$.
\endproclaim

\demo{Proof} Set $M=N^\perp$. Write $n=\mu(N)\geq \mu(M)$
and let $U$ be
a submodule of $M$ generated by $n-1$ general linear combinations
of generators of $M$. By Remark 2.2,  there exists an element $x\in
N$ such that
$U:_RM = N^*(x)$. The latter ideal has height at most $r$ by assumption.
Applying Proposition 3.7bis with $e = n-r$ and $t = n-1$ we  conclude that
$\ell(M)\geq n$, hence $\ell(M) = n$.  \qed  \enddemo

Combining Propositions 3.6 and 3.9bis one obtains the following.  Assume
that $R$ is an equidimensional universally catenary Noetherian local ring
of dimension $d > 0$ with infinite residue field, and let $N$ be a finitely
generated $R$-module that is free of constant $r$ locally on the punctured
spectrum.  For every $x \in N$ one has $\hgt(N^*(x)) \le r < d$ if and only
if $\ell(N^{\perp}) = \mu(N)$.  To see this also notice that
$\ell(N^{\perp}) \le \dim(R) + \text{rank}(N^{\perp}) - 1 = \mu(N) + d-r-1$.

\medskip

Proposition 3.9 gives a systematic way of constructing
modules of rank $r$ with elements whose order ideals have height
exceeding $r$.

\medskip

\noindent {\bf Example 3.10.}
Let $R$ be a Noetherian local ring  with infinite residue field and let $I$
be an $R$-ideal of positive height. Suppose that $I$ satisfies
$G_\infty$ and $\ell(I)\neq \mu(I)$.  Let $N=I^\perp$ and
notice that $I=N^\perp$.  By Proposition 3.9, $N$ contains an element $x$
such that the height of $N^*(x)$ is
strictly greater than the rank of $N$.  Of course $N^*(x)$ is
proper, since $N$ has no free summands.

If $I$ is the defining ideal of a monomial curve in ${\bold P}^3$
then the analytic spread of $I$ is
at most 3
by Giminez, Morales and Simis [1993].
Such an ideal satisfies $G_4$ by Herzog
[1970]. If in addition $\mu(I)=4$, then $I$ is $G_\infty$ and
$\ell(I)\neq \mu(I)$.  The module $N=I^\perp$ has rank 3, and --- at
least if the ground field is infinite --- will
have an order ideal of height 4.

To be explicit, let $I\subset R=k[a,b,c,d]$ be the defining ideal of the
monomial curve $t\mapsto (1,t^{\alpha-1}, t^{\alpha+1}, t^{2\alpha})$,
for even numbers  $\alpha> 0$. It is
easy to check that this curve lies on the smooth quadric
$ad-bc=0$, and has divisor class
$(\alpha-1, \alpha+1)$. Its ideal is thus minimally generated by $4$ elements
and all the conditions above are satisfied. (Actually the
same is true for any curve in this divisor class, monomial
or not.) The module $N$ may be explicitly described as
the image of the right-hand map in the left exact sequence
$$
0\ra
R
\overset{\left( \matrix
ad-bc \\
c^{\alpha+1}-b^{\alpha-1}d^2 \\
ac^\alpha-b^\alpha d \\
b^{\alpha +1}-a^2c^{\alpha -1}
\endmatrix \right)}\to\lra
R^4
\overset{\left( \matrix
-b^\alpha &0&a&c \\
-ac^{\alpha -1} & 0 & b &d \\
-b^{\alpha -1}d & a & -c & 0\\
-c^\alpha &b&-d&0
\endmatrix \right)}\to\lra
R^4
$$
obtained by dualizing the first two steps of the minimal
free resolution of $R/I$. From this we see at once that
the third generator of $N$ has order ideal $(a,b,c,d)$,
of height 4.
\medskip

We now give an example of a graded module $N$ in which all
homogeneous elements have order ideals of  height at most rank $N$,
although there are
inhomogeneous elements whose order ideals have bigger  height.

\medskip

\noindent{\bf Example 3.11.}  Let $k$ be an infinite field and let
$R=k[z_1,\dots,z_6]$ be a
polynomial ring, graded with all the variables in degree 1.
Set
$$
I= I_2 \pmatrix
z_1 & z_2^2& z_3^2 & 0\\
0   & z_4^2& z_5^2 & z_6^2
\endpmatrix.
$$
The perpendicular module $N=I^\perp$ has rank 5.
We claim first that all the order ideals of
homogeneous elements of $N$
have height at most 5.
By  Remark 2.2 and Corollary 2.3 it suffices  to show that
the ideal $I$ has no reductions generated by
5 homogeneous elements.

Suppose on the contrary that $J$ is a reduction of $I$
generated by 5 homogeneous elements.
The lowest degree part of $I$ is generated by 3 analytically independent
elements. The ideal $J$ must contain all three since the ideal
generated by the lowest degree elements of $J$ is a reduction
of the ideal generated by the lowest degree elements of $I$.
On the other hand, $(z_1)$ contains the lowest degree part of $I$.
Let $^{-}$ denote images in $R/(z_1)$.
The ideal $\overline J$ is
a reduction of $\overline I $. Thus $\overline I $
would have a
reduction generated by 2 elements.  However, this is impossible:
$\overline I $ is generically a complete intersection of height 2, and not
a complete intersection, so by
Cowsik and Nori [1976] any reduction of $\overline I $
has at least 3 generators.

On the other hand, if we regrade the ring
with degree$(z_1)=2$, all the generators of $I$
become homogeneous of the same degree. Since the
generators of $I$ satisfy the Pl\"ucker relation,
$I$ has a reduction $J$ generated by 5 elements homogeneous
in the new grading.  As $I$ is a complete intersection on the
punctured spectrum, the ideal $J:I$ has height 6, and thus
$N$ has an element whose order ideal has height 6 by Remark 2.2.

\medskip

We finish this section with two classes of examples arising from the Koszul
complex and the Buchsbaum-Rim complex, respectively.

Let $R=k[z_1,\dots,z_d]$ be a polynomial ring in $d>1$
variables over a field $k$, graded with the variables in degree 1,
and write $\Omega^i$ for
the $i$th syzygy module of the maximal ideal
$\m=(z_1,\dots,z_d)$.
The module $\Omega^i$ has minimal free presentation
$\wedge^{i+2}R^d \to \wedge^{i+1}R^d$; in particular
$\Omega^i/\m\Omega^i = \wedge^{i+1}k^d$.
Using the self-duality of the Koszul complex, we see at once
that $(\Omega^i)^\perp=\Omega^{d-i-2}$ for $0\leq i\leq d-2$.
Note that $\Omega^{d-1}$ is a free module.

\proclaim {Proposition 3.12} If $2\leq i\leq d-2$, or if $i=1$ and $d$ is
odd,
then all elements of $\Omega^i$ have
order ideals of height at most the rank ${d-1\choose i}$ of $\Omega^i$.
\endproclaim

\demo{Proof} We may suppose $R$ is local.
If $2\leq i\leq d-3$ there is nothing to prove, since $\Omega^i$
has no free summand and its rank is greater than the
the dimension of $R$. If $i=d-2$ then $(\Omega^i)^\perp =\m$,
which is generated by analytically
independent elements. Similarly, if $i=1$ and $d$ is odd, then
$(\Omega^1)^\perp=\Omega^{d-3}$
is generated by analytically independent elements
 by Simis, Ulrich and Vasconcelos [1993, 3.1].  Now Proposition 3.6 yields
the desired inequality in either case.
\qed
\enddemo

We can prove a more precise result for homogeneous
generators: We say that an element of $\Omega^1$ has rank $t$ if its
image in $\Omega^1/\m\Omega^1=\wedge^2k^d$ represents a linear
transformation $(k^d)^*\to k^d$ of rank $t$.  Since these linear
transformations are alternating, the rank $t$ is an even
number. Any homogeneous minimal generator
of rank $2s$ can be written as $e_1\wedge e_2+\dots+e_{2s-1}\wedge e_{2s}$,
where the $e_i$ are homogeneous minimal generators of $R^d$.

\proclaim {Proposition 3.13} The height of the
order ideal of a homogeneous generator $x$ of $\Omega^1$ is
equal to the rank of $x$. In particular, if $d$ is even, there are
elements of $\Omega^1$ with order ideals of height $d> \rank(\Omega^1) = d-1$.
\endproclaim

\demo{Proof} Let $e_1,\dots,e_d$ be homogeneous generators of $R^d$.
If $d=2s$ is even, then the generator
$e_1\wedge e_2+\dots+e_{2s-1}\wedge e_{2s}$ has rank $2s$, so the
second statement follows from the first.

The module $\Omega^1$ is the image of the map
$\wedge^2 R^d\to \wedge^1R^d$ in the Koszul complex of
$z_1,\dots, z_n$.
To prove the first statement, it suffices to consider the
order ideal of the element $x$
that is the image of $e_1\wedge e_2+\dots+e_{2s-1}\wedge e_{2s}$.
The column corresponding
to $e_{2i-1}\wedge e_{2i}$ has $\pm z_{2i}$ and $\pm z_{2i-1}$
in the $2i-1$ and $2i$ places, respectively. Thus $x$ is
mapped to an element of $R^d$ whose nonzero coordinates are
$\pm z_1,\dots,\pm z_{2s}$. Since the dual of the Koszul
complex is exact, the components of
the inclusion map $\Omega^1\to R^d$ generate all the maps
from $\Omega^1$ to $R$, and we see that the order
ideal of $x$ is $(z_1, \dots, z_{2s})$.
\qed
\enddemo

In contrast to Proposition 3.13 the next example shows that the kernel $N$
of a generic map $R^s \to R^t$ has only order ideals of height at most
$\rank(N)$ as long as $t > 1$.

\proclaim{Proposition 3.14} Let $(R,\m)$ be a local Gorenstein
ring and let $t\leq s$ be integers. Let $\chi$ be a $t$ by $s$ matrix with
entries in $\m$ and $\hgt (I_t(\chi)) =s-t+1$. Set $N= ker(\chi)$.
Except in the case where $t=1$ and $s$ is even,
$\hgt(N^*(x))\leq \rank(N)$ for every $x\in N$.
\endproclaim

\demo{Proof} We may assume that $s-t = \rank(N) > 1$.
Let $x_1,\dots,x_n$ be a minimal generating set of $N$, let
$S=R[Z_1,\dots,Z_n]$ be a polynomial ring, and write
$y=\sum_{i=1}^nZ_ix_i\in N\otimes_RS$. Further, let
$I=(N\otimes_RS)^*(y)$ be the order ideal of $y$
and $J$ its unmixed part.
By Proposition 2.4(2) it suffices to show that $\hgt(I_{\m S})\leq \rank(N)$.

To this end write $M=N^\perp$. As in the proof of Remark 2.2 one sees
that $\Sym(M)\cong S/I$. Since $M$ satisfies $G_\infty$, one
has $\dim(\Sym(M))= \dim(R) +\rank(M)$ by Huneke and Rossi [1986, 2.6],
which gives $\hgt(I)=\rank(N)=s-t$. In this setting,
Migliore, Nagel and Peterson show in [1999, 1.5b] that $I=J$ if
$s-t$ is even, whereas
$J/I\cong \Sym_{(s-t-1)/2}$(coker$(\chi))\otimes_RS$ if
$s-t$ is odd. Thus for $s-t > 1$ odd and $t > 1$ one has $\mu_{\m S}(J/I) >
1$, which yields $J_{\m S}\neq S_{\m S}$. Hence in either case
$$
\hgt(I_{\m S})=\hgt(J_{\m S})=\hgt(J)=\hgt(I)=\rank(N).\qed
$$
\enddemo

\bigskip
\goodbreak
\centerline{\bf 4. Order Ideals of Low Rank Modules\rm}
\bigskip

It turns out that order ideals of elements in modules of low rank
are particularly well-behaved. For example if  $N$ is a module of
rank 1, then (modulo torsion) $N$ is isomorphic to an ideal $I$
containing a nonzerodivisor.  If $x\in I$ is a nonzerodivisor of
$R$ then $I^*(x)=(x):_RI$ which is either the unit ideal or of
grade  1. If on the other hand $x$ is a zerodivisor contained in
an associated prime $Q$ of $R$, then $I^*(x)$ is also contained in
$Q$ and thus has grade 0. The following propositions extend this
kind of result to modules of rank 2 and 3 as well as $k^{th}$
syzygies of rank $k$ having finite projective dimension. The case
of rank 2 modules over regular local rings had already been
treated in Evans and Griffith [1982, p.377].

\proclaim{Proposition 4.1} Let $R$ be a Noetherian ring and let
$N$ be a finitely generated $R$-module.  Assume either
\roster
\item $N$ is orientable of rank 2; or
\item  $R$ contains a field and $N$ is a $k^\th$
syzygy of rank $k$ having finite projective dimension.
\endroster
If $x\in N$ then  $N^*(x)=R$ or grade$(N^*(x))\leq \rank(N)$.
\endproclaim

\demo{Proof} We may assume that $R$ is local.
Consider the exact sequence
$$
0\ra Rx\ra N\ra X\ra 0.
$$
We suppose that grade$(N^*(x))>\rank(N)$.
It follows that the annihilator of $x$ is 0, and
we will show that $X$ is free. Hence $N^*(x)=R$ as
required.

In case (1) we may assume that $N = N^{**}$.
If we localize the above sequence at an arbitrary prime $Q$
with depth$(R_Q)\leq 2$, then the sequence splits and hence $X_Q$ is
reflexive.
If depth$(R_Q)\geq 3$ then
depth$(X_Q)\geq 2$. It follows that  $X$ is orientable and reflexive of
rank one, hence free, completing the proof in this case.

Now suppose we are in case (2).
After localizing at a prime $Q$ with depth$(R_Q) \leq k$,
the above sequence splits.
Since $N_Q$ is free by the Auslander-Buchsbaum formula, it follows
that $X_Q$ is free as well for any such $Q$.  If $Q$ is
any prime with depth$(R_Q) >  k$, then depth$(N_Q)\geq k$ and
therefore depth$(X_Q)\geq k$.    Thus by
Hochster and Huneke [1990, 10.9], the module $X$ is a $k$th syzygy.
It has rank $k-1$ and finite projective dimension, so the
version of the Evans-Griffith Syzygy Theorem due
to Hochster and Huneke [1990, 10.8] and Evans and Griffith [1989, 2.4]
 implies that $X$ is free.
\qed
\enddemo

The assumption of orientability in Proposition 4.1 is necessary:
Let $k$ be a field, $R=k[Z_0,\dots,Z_3]/(Z_0Z_3-Z_1Z_2)$, and let
$z_i$ denote the image of $Z_i$ in $R$.
Let $M$ be the ideal $(z_0^2, z_0z_1,z_1^2)$.
If $N=[z_0^2, z_0z_1,z_1^2]^\perp$ then $N$ is an
$R$-module of rank 2. By Remark 2.2
the element $(z_0z_1)^\perp$ has order ideal
$(z_0,\dots,z_3)$, which has height 3.

\proclaim{Proposition 4.2} Let $R$ be a Gorenstein ring and let $N$ be an
orientable $R$-module of rank $3$ that satisfies $S_3$ and is free in
codimension $2$. If $x\in N$ then either $N^*(x) = R$ or
$\hgt (N^*(x))\leq 3$.
\endproclaim

\demo{Proof} We may assume that $R$ is local, and we write $n = \mu(N)$, $e
= n-3$. By Proposition 4.1 we may
assume that $N$ has no nontrivial free summand. We will prove that
$\hgt (N^*(x))\leq 3$.

Suppose the contrary and write $M = N^{\bot}$. By Theorem 3.1, $x$ can be
extended to a minimal generating set of $N$.
 Using Remark 2.2, a row ideal in some minimal presentation matrix
of $M$ then has height $> 3$. Let $u_1, \dots ,u_n$ be generic elements in
$M\otimes_R{R'}$ defined over a local ring $R'$ that is obtained from $R$ by a
purely transcendental residue field extension, and set $F =
R'u_1+\cdots+R'u_{e-1} \subset U = R'u_1+\cdots+R'u_{n-1} $. The genericity
of $u_1, \dots ,u_n$ implies that
$\hgt (U:_{R'}(M\otimes_R{R'})) > 3$ as shown in Proposition 2.4(1).
To simplify notation
we will write $R = R'$.

Because $N$ is $S_3$, the Acyclicity Lemma of Peskine and Szpiro [1972]
shows that  \linebreak $\ext^1(N^*,R) = 0$. From the exact sequence
(*) of Proposition 2.5 we see that  $\ext^2(M,R) = 0$.
Because $N$ has no free summands, $\mu(M) = n$ and hence $U \ne M$.
The module $M$ is free locally in
codimension $2$, and by Proposition 2.5 it is orientable of rank $e$.
Thus by Simis, Ulrich and Vasconcelos [1998, 3.2], $F$ is free, and $M/F$
is isomorphic to an ideal $I\subset R$, with $\hgt(I) \geq 2$,
that is a complete intersection locally
in codimension two. Clearly $\ext^3(R/I,R)\cong \ext^2(I,R)\cong
\ext^2(M,R) = 0$.
Let $J$ be the image of $U$ in $I$. Notice that $\mu(J)\leq (n-1)-(e-1) = 3$.
Furthermore $J\ne I$ and $\hgt (J:I) > 3$, because $I/J\cong M/U$.

If $\hgt (I) \geq 3$ then $J$ is a complete intersection and $\hgt
(J:I)\leq 3$, and we are done.
Otherwise $\hgt (I) = 2$. Then the condition  $\ext^3(R/I,R) = 0$
implies that the factor ring of $R$ by any link of $I$ satisfies $S_2$;
see Chardin, Eisenbud and Ulrich [1998, 4.4].  Since we are in a
case where
$I$ has height 2 and is a complete intersection in codimension 2,
and $\mu(J)\leq 3$,
we can  apply Chardin, Eisenbud and Ulrich
[1998, 3.4]
to obtain $\hgt (J:I) \leq 3$ as required.
\qed
\enddemo

In Proposition 4.2 the assumption of freeness in codimension 2 can be
weakened to requiring   that $\rf_Q(N) \ge 2$ whenever $\dim(R_Q) = 2$.
However, the $S_3$ condition is necessary,
as can be seen from  the monomial curves discussed in Example 3.10.

\medskip

We can use Proposition 4.1 to prove that, under a vanishing hypothesis
on some $\ext^i(M,R)$, the colon ideal $U:_RM$ has at most the
\it expected \rm grade. In preparation, recall that if
$M^*$ has a rank then $M$ does too.

\proclaim{Proposition 4.3} Let $R$ be a Noetherian local ring containing
a field,  let $M$ be a finitely generated torsion free $R$-module
such that $M^*$ and $\ext^1(M,R)$ have
finite projective dimension,
and set $e = \text{rank}(M)$.
If $s$ is an integer such that
$\ext^{i}(M,R) = 0$ for $2\leq i\leq s-1$
then, for every submodule $U\subsetneq M$ with $\mu(U)-e+1\leq s$,
$$
\text{grade}(U:_RM)\leq \max\{0,\mu(U)-e+1\}.
$$
\endproclaim

\demo{Proof}   We may assume that rank$(U) = e$ since otherwise
$\gr(U:_RM) = 0$. Then $s\geq \mu(U)-e+1\geq 1$. Lowering $s$ if necessary,
it suffices to prove that
$\gr(U:_RM)\leq s$.
If $s = 1$, then $U$ is
free and thus $\gr(U:_RM)\leq 1 = s$. Therefore we may assume that
$s\geq 2$. Suppose that $\gr(U:_RM) > s$.
Choose a submodule $U\subsetneq V \subset M$ such that $V/U$ is cyclic.
Clearly
$\gr(U:_RV) > s$ and $\gr(V:_RM) > s$.
The latter inequality
implies that
$\ext^i(V,R)\cong \ext^i(M,R)$ for every $0\leq i\leq s-1$, forcing $V^*$ and
$\ext^1(V,R)$ to have finite projective dimension and
$\ext^i(V,R) = 0$ for $2\leq i\leq s-1$. We are free to replace $M$ by $V$
and from now on we assume that $M/U$ is cyclic.

Since $M/U$ is cyclic and $\mu(U)\leq s+e-1$, there exists a generating set
$u_1, \dots ,u_n$ of
$M$ such that the first $n-1$ elements generate $U$ and $n = s+e$. The
module $M$ is torsionless since it is torsion free and has a rank. By the
remark at the
end of Section 2, there exists a finitely generated module $N$ with
generators $x_1, \dots ,x_n$ so that $M = [x_1,\dots ,x_n]^{\bot}$ and
$u_i = (x_i)^{\bot}$.  By Remark 2.2, $U:_RM = N^*(x_n)$.
Using the sequence (*) of Proposition 2.5 we see that $N^{**}$ has finite
projective dimension since $M^*$  and  $\ext^1 (M,R)$ have finite
projective dimension,  and that $N^{**}$ is an $s^\th$ syzygy since
$\ext^i(M,R) = 0$ for $2\leq i\leq s-1$.
As rank$(N^{**}) = n-e = s$, Proposition 4.1(2) now proves that
$\text{grade}(N^*(x_n))\leq s$, giving the conclusion.
\qed\enddemo

\proclaim{Corollary 4.4} Let $R$ be a regular local ring containing a
field and let $I$ be an ideal satisfying $\ext^i(R/I,R) = 0 $ for
$3\leq i\leq s$. Then for every ideal $J\subsetneq I$ with $\mu(J)\leq s$,
$\hgt(J:I)\leq \mu(J)$.\qed
\endproclaim

Note that Corollary 4.4 is interesting only if the
height of $I$ is one or two, and reduces to the height two case.
One should compare the results of Chardin, Eisenbud, and
Ulrich [1998, 3.4 and
4.2], which yield the same conclusion for ideals of any
height $g$ under the (incomparable) assumptions
that $I$ satisfies $G_{s}$
and  $\ext^i(R/I^{i-g},R) = 0$ for $g+1 \leq i \leq s$.
We did not expect a result that  avoids  reference
to the powers of I!

\proclaim{Corollary 4.5} Let $R$ be a Noetherian local ring
containing a field, let $M$ be a finitely generated torsion free $R$-module
such that $M^*$ and $\ext^1(M,R)$ have finite projective dimension,
and set rank$(M) = e$.
If $s$ is an integer such that $R$ satisfies $S_{s+1}$,
the module $M$ satisfies
$G_{s+1}$,
and  $\ext^i(M,R) = 0$ for $2\leq i\leq s-1$,
then
$$
\ell(M)\geq \text{min}\{\mu(M), e+s\}.
$$
\endproclaim

\demo{Proof} We may assume that $R$ has an infinite residue field. Write $n
= \mu(M)$ and $t = \text{min}\{n,e+s\} - 1$. We may suppose that $t \ge 1$.
Let
$U$ be a submodule of $M$ generated by $t$ general linear combinations
of generators of $M$. As $t\leq n-1$ one has  $U\ne M$ and therefore
$\hgt(U:_RM)\leq t-e+1$ by Proposition 4.3. But then Proposition 3.7bis
implies that $\ell(M)\geq t+1$.
\qed\enddemo

\goodbreak
\bigskip\centerline{\bf 5. Heights of Trace Ideals\rm}
\bigskip

The surprising fact pursued in this section may be
informally summarized by saying that if an order
ideal of an element in a module $N$
is ``bigger than it should be'', then  the trace ideal of $N$
is not much larger than this order ideal.

\proclaim{Proposition 5.1} Let $R$ be a Noetherian ring, let  $N = Ry_1 +
\cdots + Ry_n$ be an $R$-module, and let $x,y$ be elements of $N$. Write
$Y_i = N/Ry_i$ and $Y = N/Ry$.
\roster
\item If $\hgt(N^*(x)/Y^*(x)) > 1$, then $N^* (y) \subset \sqrt{N^*(x)}$.
\item If $\hgt(N^*(x)/Y^*_i(x)) > 1$ for $1 \le i \le n$, then
$\sqrt{\text{tr}(N)} = \sqrt{N^*(x)}$.  In particular $\hgt(\text{tr}(N))
\le \mu(N^*), $ unless  $\text{tr}(N) = R$.
\endroster
\endproclaim

\demo{Proof} To prove (1) suppose there exists a prime ideal $Q$ of $R$
with $N^*(x) \subset Q$, but $N^* (y) \not\subset Q$.  Replacing $R$ by
$R_Q$ we may then assume that $N^* (x) \neq R$ and $N = R \oplus Y$.
Writing $x = r + z$ with $r \in R, z \in Y,$   we obtain $R \neq N^*(x) =
Rr + N^*(y)$.  Thus $\hgt(N^*(x)/N^*(y)) \le 1$, which yields a
contradiction.  This proves (1).

Part (2) is an immediate consequence of (1) since $N^*(x) \subset
\text{tr}(N) = N^*(y_1) + \cdots + N^*(y_n)$.
\qed
\enddemo

\remark{\bf Remark 5.2} The height assumptions in Proposition 5.1 are
automatically satisfied if $R$ is locally equidimensional and catenary and if
$\hgt(N^*(x)) > \bight(Y^*(x)) + 1$  (for (1)) or   $\hgt(N^*(x)) > \bight
(Y^*_i(x)) + 1$ (for (2)).
\endremark

\medskip

There is a corresponding statement for colon ideals.  When combined with
Remark 2.2 it could be used to deduce Proposition 5.1:

\proclaim{Proposition 5.3}  Let $R$ be a Noetherian ring and let $M$ be a
finitely generated $R$-module.  Let $U = Ru_1 + \cdots + R u_s \subset M$
be a submodule and write $U_i = Ru_1 + \cdots + Ru_i$.  If $\hgt(U_i:_R
M/U_{i-1}:_R M) > 1$ for $1 \le i \le s$, then $\sqrt{\text{Fitt}_{j+s}(M)}
= \sqrt{\text{Fitt}_j(M/U)}$ for any $j \ge 0$.  In particular $V :_R M
\subset \sqrt{U :_R M}$ for every $s$-generated submodule $V \subset M$.
\endproclaim

\demo{Proof} The asserted equality is equivalent to the statement that
$u_1, \dots, u_s$ form part of a minimal generating set of $M$ locally at
each prime in the support of $M/U$.  Thus it suffices to prove that for $1
\le i \le s$, the image of $u_i$ is a minimal generator of $M/U_{i-1}$
locally on the support of $M/U_i$.  This reduces us to the case $s = 1$.

Suppose $u=u_1$ is not a minimal generator of $M$ locally at a prime $Q$ in
the support of $M/U = M/U_1$.  Replacing $R$ by $R_Q$ we may assume
$(R,\m)$ is local, $M \neq 0$, and $u \in \m M$. If $\var$ is a matrix with
$n$ rows presenting $M$, then $I_n(\var) \neq R$ and we obtain a
presentation matrix $\psi$ of $M/U$ by adding one column with entries in
$\m$.  Thus by Bruns [1981, Corollary 1]
(see also Eisenbud and Evans [1976, 2.1]),
$\hgt(I_n(\psi)/I_n(\var))
\leq 1$, which gives $\hgt(\text{ann}(M/U)/\text{ann}(M)) \leq 1$, contrary
to our assumption.
\qed
\enddemo

In  general one has the inclusion $\text{Fitt}_j(M/U) \subset
\text{Fitt}_{j+s} (M)$ for any $s$-generated submodule $U$ of a finitely
generated module $M$.  One may ask which power of $\text{Fitt}_{j+s}(M)$ is
contained in $\text{Fitt}_j(M/U)$ under the assumptions of Proposition 5.3.

The height assumption in Proposition 5.3 implies that
$\hgt(U :_R M/0:_R M) \ge 2s$.
On the other hand a lower bound for $\hgt(U :_R M/0 :_R
M)$ alone does not suffice to deduce the equality
$\sqrt{\text{Fitt}_{j+s}(M)} = \sqrt{\text{Fitt}_j(M/U)}$.  For instance,
let $R = k[z_0, \dots, z_n]$ be a polynomial ring over a field, $\var$ the
$3$ by $n$ matrix
$$\left( \matrix
0 & 0 & \cdots & 0 \\
z_0 & 0 & \cdots &0 \\
z_1 & z_2 & \cdots & z_n
\endmatrix \right),$$
 which gives a map from $R^n$ to $R^3 = Re_1 \oplus Re_2 \oplus Re_3$, $M$
the cokernel of $\var$, and \linebreak $U \subset M$ the submodule
generated by the images of $e_1$ and $e_2$.  In this case, one has
ht$(U :_R M/0 :_R M) = n$, whereas $\text{Fitt}_2(M) = (z_0,
\dots, z_n) \not\subset \sqrt{\text{Fitt}_0(M/U)} = (z_1, \dots,
z_n)$.

\proclaim{Proposition 5.4} Let $R$ be a universally catenary Noetherian
ring and let $N$
be a finitely generated $R$-module. Assume that one of the following
conditions hold:
\roster
\item $R$ satisfies $S_3$ and $N$ is orientable
of rank 3; or
\item $R$ is Gorenstein and $N$ is  orientable of rank $4$, satisfies
$S_3$ and is free in codimension $2$; or
\item $R$ satisfies $S_{k+1}$ and contains a field, and $N$ is a $k$th
syzygy of rank $k+1$ having finite projective dimension.
\endroster
If $x\in N$ satisfies $\hgt(N^*(x)) > \text{rank}(N)$ then
$\rad{\text{tr}(N)} = \rad{N^*(x)}$.
\endproclaim

\demo{Proof} We may assume that $R$ is local,
$N$ is torsion free and nonzero, and $N^*(x)\neq R$.
If we are in case (3) with $k=0$ then
$N$ is orientable of rank 1, so $N^*\cong R$
and $N^*(x)$ has height 1, contradicting the hypothesis.
In all other cases we may suppose
that $R$ is $S_2$, hence equidimensional.

After passing to a purely transcendental extension of the
residue field of $R$, we consider generic generators $y_1, \dots ,y_n$ of $N$
and set $Y_i = N/Ry_i$. By Proposition 5.1(2) and Remark 5.2,
we only need to prove that
$\bight(Y_i^*(x))\leq \text{rank}(N) - 1$ for $1\leq i\leq n$. Write
$y = y_i, Y = Y_i$, and notice that $Y^*(x)\subset N^*(x)\ne R$. By
Proposition 2.4(2)
$\hgt(N^*(y))\geq \hgt(N^*(x)) > \text{rank}(N)$.
Thus $Ry \cong R$ and assumptions (1), (2), and (3) pass from
$N$ to $Y$, except that $\text{rank}(Y) = \text{rank}(N)-1$. After localizing
at minimal primes of $Y^*(x)$ we may apply Propositions 4.1 and 4.2 to
conclude that $\bight(Y^*(x))\leq \text{rank}(Y) = \text{rank}(N)-1$.
\qed
\enddemo

We now consider bounds for
the height of the trace ideal $\tr(N)=N^*(N)$. Of course
if $R$ is regular local and
the height of the order ideal $N^*(x)$ is bounded by $\rank(N)$
for every $x$ then $\hgt(\tr(N))\leq ri$ because the
heights of ideals are subadditive in a regular local ring
(Serre [1958, V.6.3]). In general, however, no such
strong inequality is true. For the module $N$ of rank $r$ which is
the image of the generic $m\times n$ matrix of rank $r$, all the
heights of order ideals are bounded by $r$ but
the height of the trace ideal is $r(m+n-r)$ (this is
the height of the $1\times 1$ minors of the generic matrix
minus the height of the $r+1)\times r+1$ minors of that
matrix---see Remark 2.2.  Bruns [1981, Corollary 1]
shows that this is a universal bound.

Curiously, if we assume
that some element has an order ideal larger than the rank
(which would seem to make the trace ideal larger) then we
can get a bound which is asymptotically much sharper (as
$(m+n)/r\to\infty$ it has order $(r-2)m+(r-3)n$ instead of
$r(m+n)$).
It is convenient at the same
time to give a bound
on $N^*(U)$ for any $i$-generated submodule $U$:

\proclaim{Theorem 5.5} Let $R$ be an equidimensional universally
catenary Noetherian local ring and let $N$ be an orientable $R$-module of
rank $r$.  Let $U \subset N$ be a
submodule not containing any nontrivial free summand of
$N$, and write $i = \mu(U), ~ m = \mu(N^*)$.  Then either
$\hgt(N^*(x)) \leq r$ for every $x \in U$ or else

\roster
\item $\hgt(N^*(U)) \leq (r-2)(m-r+3) + (r-3)i$, in case $R$ satisfies $S_3$;
\item $\hgt(N^*(U)) \leq (r-3)(m-r+3) + (r-3)i$, in case $R$ is Gorenstein,
and $N$ satisfies $S_3$ and is free in codimension 2;
\item $\hgt(N^*(U)) \leq (r-k)(m-r+3) + (r-3)i$, in case $R$ contains a
field and satisfies $S_{k+1},$ and $N$ is a $k$th syzygy of finite
projective dimension.
\endroster
\endproclaim

\demo{Proof}  Assume that $\hgt(N^*(x))> r$
for some $x\in U$. We will prove (1) and (2) by induction on $i$.  Using
Propositions 4.1 and 4.2 we see that $r \geq 3$ in (1),
$r\geq 4$ in (2), and $r \geq k+1$ in (3).

By Theorem 3.1 we can write $U = Rx_1+ \cdots +Rx_i$ with $\hgt(N^*(x_1)) >
r$.
If $i = 1$, then $\hgt(N^*(U)) = \hgt(N^*(x_1))\leq m$, which yields the
desired estimates in this case (as $r \geq 3$, or $r\geq 4$, or $r \geq k+1$,
respectively).  Hence we may assume $i\geq 2$. Write $V = Rx_1+\cdots
+Rx_{i-1}$.
If $N^*(x_i)\subset \rad{N^*(V)}$ we can replace $U$ by $V$ and apply the
induction hypothesis  (note that $x_1\in V$ and $\hgt(N^*(x_1)) > r$).
Otherwise, we may choose a prime $Q$ containing
$N^*(V)$ such that $N^*(x_i)\not\subset Q$.

Set $X = N/Rx_i$, let $z_j$ be the image of $x_j$ in $X$ for
$1\leq j\leq i-1$, and  let $Z$ be the image of $V$ in $X$. We estimate
the height of $N^*(U)$ as follows:
$$\hgt(N^*(U)) = \hgt(N^*(x_i) + N^*(V))\leq m + \hgt(N^*(V))\leq m +
\hgt(N^*(V)_Q).$$
Since $N^*(x_i)\not \subset Q$, $N_Q \cong R_Q \oplus X_Q$. We may write
$x_j = r_j + z_j$ under this isomorphism. We deduce $N^*(x_j)_Q = (r_j,
X^*(z_j))_Q$.  It follows that
$$
N^*(V)_Q = (r_1,\dots ,r_{i-1})_Q + X^*(Z)_Q
$$
and $\hgt(X^*(z_1)_Q)\geq \hgt(N^*(x_1)_Q) - 1 > r - 1$.
Furthermore, as $N^*(V)_Q\ne R_Q$
we may estimate the height of $N^*(V)_Q$; combining with the
estimate above we obtain,
$$\hgt(N^*(U))\leq m + (i-1) + \hgt(X^*(Z)_Q).$$
We apply our induction hypothesis to $Z_Q\subset X_Q$ noting that
$\mu(Z_Q)\leq i-1$, $X_Q$ is orientable of
rank $r-1$, $\mu(X_Q^*) \leq m-1$, and $\hgt(X_Q^*(z_1)) > r-1$.

We formally set $k=2$ for (1), $k=3$ for (2), and $k=k$ for (3).
By induction,
$$
\hgt(N^*(U)) \leq m + (i-1) + ((r-1)-k)((m-1)-(r-1)+3) + ((r-1)-3)(i-1)
$$
and the desired formulas follow.
\qed
\enddemo

For example, let $I$ to be the ideal of the curve $t\mapsto
(1,t,t^3,t^4)$ in ${\bold P}^3$ treated in Example 3.10. If
$N=I^\perp$, then the inequality of Proposition 5.1(1) is sharp for
every $i$ (here $r=3, m=n=4$).
We do not have examples of rank $\geq 4$ where the inequality is sharp.

\vfill\eject

\centerline{\bf Bibliography}
\bigskip
\refstyle{A}

\enddocument